\documentclass[12pt]{article}

\usepackage{amssymb}
\usepackage{amsmath}
\usepackage{graphicx}
\usepackage{hyperref}

\numberwithin{equation}{section}

\begin{document}

\newcommand{\bd}{\begin{displaymath}}
\newcommand{\ed}{\end{displaymath}}
\newcommand{\ds}{\displaystyle}
\newcommand{\bp}{\underline{\bf Proof}:\ }
\newcommand{\ep}{{\hfill $\Box$}\\ }
\newcommand{\be}{\begin{equation}}
\newcommand{\ee}{\end{equation}}
\newcommand{\ba}{\begin{array}}
\newcommand{\ea}{\end{array}}
\newcommand{\bea}{\begin{eqnarray}}
\newcommand{\eea}{\end{eqnarray}}
\newcommand{\nt}{\noindent}

\newtheorem{0}{DEFINITION}[section]
\newtheorem{1}{LEMMA}[section]
\newtheorem{2}{THEOREM}[section]
\newtheorem{3}{COROLLARY}[section]
\newtheorem{4}{PROPOSITION}[section]
\newtheorem{5}{REMARK}[section]
\newtheorem{6}{EXAMPLE}[section]
\newtheorem{7}{ALGORITHM}[section]
\newtheorem{8}{CONJECTURE}[section]

\title{An Alternative Framework for Irreducibility and Primitivity of Nonnegative Tensors}
\author{
Jianhong Xu\thanks{School of Mathematical and Statistical Sciences, Southern Illinois University Carbondale, Carbondale, IL 62901, USA. Email: \texttt{jhxu@siu.edu}} 
% \and  \thanks{}
% \and  \thanks{}
}

\maketitle

\begin{abstract}  
Motivated by some recent studies on higher order Markov chains and well-known characterizations for irreducibility and primitivity of nonnegative matrices, we propose in this paper an alternative framework for irreducibility and primitivity of nonnegative tensors, giving rise to two new notions: s-irreducibility and s-primitivity. Our framework includes the relevant results on nonnegative matrices as its special cases, yet it expands existing results regarding irreducibility and primitivity of nonnegative tensors. In addition to its tensor theoretic significance, such a framework has important implications for applied problems such as higher order Markov chains.
\end{abstract}

\nt {\bf Keywords}: 
nonnegative tensors, tensor product, irreducibility, primitivity, accessibility, aperiodicity

\nt {\bf AMS Subject Classification}: 15A69, 15A72, 15B48, 46B28 

\section{Introduction}
\label{intro}
\setcounter{equation}{0}

Let us start by giving a few remarks about notation. Throughout this work, we shall denote scalars and vectors by lowercase letters, matrices by uppercase letters, and tensors by calligraphic uppercase letters. For the entries of a vector, matrix, or tensor, we shall denote by the matching lowercase letter with subscripts $i$, $j$, etc. for the indices. As an exception, however, we shall denote a zero vector, matrix, or tensor simply as $0$, while inequalities such as vector $x \ge 0$ and tensor ${\cal A} > 0$ will be interpreted in an entrywise sense. In addition, we shall assume that $m, n \ge 2$ and $\alpha, \beta, \gamma \ge 0$ are integers. For simplicity, we shall denote $\{1, 2, \ldots, n\}$ as $\langle n \rangle$. The $i$th columns of an $n \times n$ identity matrix $I_n$ will be denoted by $e_i$.

Irreducibility and primitivity are important notions in the study of nonnegative square matrices; see, for example, \cite{BP, HJ}. Specifically, we have: 

\begin{0}
\label{irr1}
Let $A=[a_{ij}]$ be a nonnegative $n \times n$ matrix. Denote $A^\alpha=[a^{(\alpha)}_{ij}]$. If given any $i, j \in \langle n \rangle$, there exists $\alpha \ge 1$, which may depend on $i$ and $j$, such that $a^{(\alpha)}_{ij}>0$, then $A$ is called irreducible.
\end{0}

\begin{0}
\label{pri1}
Let $A$ be a nonnegative $n \times n$ matrix. If there exists $\alpha \ge 1$ such that $A^\alpha>0$, then $A$ is called primitive.
\end{0}

From these definitions, it is clear that primitivity implies irreducibility, but the reverse implication does not hold.

In the context of the classical first order (homogeneous) Markov chains, irreducibility and primitivity are usually called ergodicity and regularity, respectively \cite{Ios, KS}. In fact, the theory of first order Markov chains has been intertwined with that of nonnegative matrices \cite{BP}. Naturally, such a connection carries over to the study of higher order (homogeneous) Markov chains and that of nonnegative tensors \cite{CZ, GLY, HQ, LZ, LN, WC}. This is a major thrust behind this work.

Definition \ref{irr1} can equivalently be given in terms of the zero-nonzero pattern of the entries of $A$ or the strong connectivity of the underlying digraph associated with $A$ \cite{HJ}. Speaking of the former, it can be stated as:

\begin{0}
\label{irr2}
Let $A=[a_{ij}]$ be a nonnegative $n \times n$ matrix. If given any $J$ such that $\emptyset \ne J \subsetneq \langle n \rangle$, $a_{ij}>0$ for some $i \in J$ and $j \in J^c$, then $A$ is called irreducible.
\end{0}

This leads to the following counterpart for nonnegative tensors \cite{Lim}.

\begin{0}
\label{irr3}
Let ${\cal A}=[a_{i_1i_2\ldots i_m}]$ be a nonnegative $m$th order, $n$ dimensional tensor. If given any $J$ satisfying $\emptyset \ne J \subsetneq \langle n \rangle$, $a_{i_1i_2\ldots i_m}>0$ for some $i_1 \in J$ and $i_2, \ldots, i_m \in J^c$, then $\cal A$ is called irreducible.
\end{0}

On the other hand, it is obvious that Definition \ref{pri1} can be rephrased as:

\begin{0}
\label{pri2}
Let $A$ be a nonnegative $n \times n$ matrix. If there exists $\alpha \ge 1$ such that $A^\alpha x>0$ for any $n$-vector $x$ satisfying $0 \ne x \ge 0$, then $A$ is called primitive.
\end{0}

This motivates a way of generalizing primitivity to nonnegative tensors. Let ${\cal A}$ be an $m$th order, $n$ dimensional tensor. For any $n$-vector $x$, $T_{\cal A}(x)$ is defined to be an $n$-vector whose entries are given by 
\be
\label{map}
[T_{\cal A}(x)]_i=\sum_{i_2, \ldots, i_m \in \langle n \rangle}a_{ii_2\ldots i_m}x_{i_2}\cdots x_{i_m}, ~i \in \langle n \rangle.
\ee
Observe that $T_{\cal A}(x) \ge 0$ whenever ${\cal A} \ge 0$ and $x \ge 0$. Moreover, $T_{\cal A}(cx)=c^{m-1}T_{\cal A}(x)$ for any scalar $c$. In the literature, $T_{\cal A}$ is often chosen as \cite{CPZ}
$$[T_{\cal A}(x)]_i=\left(\sum_{i_2, \ldots, i_m \in \langle n \rangle}a_{ii_2\ldots i_m}x_{i_2}\cdots x_{i_m}\right)^{1/(m-1)}\!\!\!\!\!, ~i \in \langle n \rangle.$$
This modification, however, does not affect the relevant results in what follows. With the map $T_{\cal A}$ in (\ref{map}), we now quote \cite{CPZ}: 

\begin{0}
\label{pri3}
Let $T_{\cal A}(x)$ be associated with a nonnegative $m$th order, $n$ dimensional tensor $\cal A$ and given by (\ref{map}). If there exists $\alpha \ge 1$ such that 
$$(\underbrace{T_{\cal A}\circ T_{\cal A} \circ \cdots \circ T_{\cal A}}_{\alpha})(x)=\underbrace{T_{\cal A}(\cdots (T_{\cal A}(T_{\cal A}}_{\alpha}(x))\cdots)>0$$
for any $n$-vector $x$ satisfying $0 \ne x \ge 0$, then $\cal A$ is called primitive.
\end{0}

For brevity, we shall denote $\underbrace{T_{\cal A}\circ T_{\cal A} \circ \cdots \circ T_{\cal A}}_{\alpha}$ as $T_{\cal A}^{\,\circ \alpha}$ in the sequel.

Analogous to the matrix case of $m=2$, if a nonnegative tensor is primitive, then it is irreducible \cite{QL}.

Definitions \ref{irr3} and \ref{pri3} concerning irreducibility and primitivity of nonnegative tensors are critical since they have been part of the foundation for investigating a variety of problems such as directed hypergraphs, the generalized Perron-Frobenius theory, $Z$-eigenvalues and $Z$-eigenvectors, and convergence of relevant numerical methods; see, for example, \cite{CPZ, CLN, QL} and the references therein. Note also that these definitions reduce to Definitions \ref{irr2} and \ref{pri2} as special cases when $m=2$.

Definitions \ref{irr1} and \ref{pri1}, on the other hand, point to an alternative framework for irreducibility and primitivity of nonnegative tensors, i.e., to specify these concepts using some well-defined, meaningful tensor power. In light of this thinking, the main impetus for this work comes from recent studies on higher order Markov chains in \cite{HX24a, HX26, Xu26a, Xu26b, Xu26c}.

In particular, the tensor ``box'' product and its related power have been introduced in \cite{HX24a,HX26, Xu26c}. Such operations have practical probabilistic interpretations. Meanwhile, they have paved the way of extending notions such as ergodicity and regularity to higher order Markov chains. Although these results are formulated for stochastic transition tensors, the same methodology is applicable to a more general setting involving nonnegative $m$th order, $n$ dimensional tensors too.

The primary goals of this work are to use the tensor ``box'' product and power to develop an alternative framework for irreducibility and primitivity of nonnegative tensors and to establish some essential properties arising in this framework. To avoid any confusion with the existing terminology in the literature, we shall name these new types of irreducibility and primitivity as ``s-irreducibility'' and ``s-primitivity'', respectively. Similar to Definitions \ref{irr3} and \ref{pri3}, s-irreducibility and s-primitivity coincide with, respectively, Definitions \ref{irr1} and \ref{pri1} for the special case of matrices at $m=2$.

Before proceeding, we mention that for second order Markov chains, i.e., third order stochastic transition tensors, results similar to Theorems \ref{komo} and \ref{irrpri} can be found in \cite{Vla84, Vla85}. These results have also served as a substantial source of inspiration for this work, and will be broadened in this work so as to deal with general nonnegative $m$th order, $n$ dimensional tensors. 

The material of this work is organized as follows. In Section \ref{main}, we shall introduce the notions of s-irreducibility and s-primitivity and shall develop several essential results regarding these properties, including their relationship with the existing irreducibility and primitivity, and a full characterization of the relationship between s-irreducibility and s-primitivity, along with a number of illustrative examples. In Section \ref{concl}, we shall give a few remarks, summarizing the results of this work and suggesting some possible topics for future work.

\section{Main Results}
\label{main}
\setcounter{equation}{0}

\begin{0}
Given $m$th order, $n$ dimensional tensors ${\cal A}=[a_{i_1i_2\ldots i_m}]$ and ${\cal B}=[b_{i_1i_2\ldots i_m}]$, the ``box'' product of $\cal A$ and $\cal B$, denoted as ${\cal A}\boxtimes {\cal B}$, is an $m$th order, $n$ dimensional tensor ${\cal C}=[c_{i_1i_2\ldots i_m}]$, whose entries are determined by 
\be
\label{bprod}
c_{i_1i_2\ldots i_m}=\sum_{j \in \langle n \rangle} a_{i_1ji_2\ldots i_{m-1}}b_{ji_2\ldots i_m}, ~i_1, i_2, \ldots, i_m \in \langle n \rangle.
\ee
\end{0}

Clearly, when $m=2$, ${\cal A}\boxtimes {\cal B}$ reduces to the usual multiplication of two $n \times n$ matrices. For the more general $m \ge 3$ case, however, a key difference arises, namely, the $\boxtimes$ product is usually not associative. In other words, for $m$th order, $n$ dimensional tensors $\cal A$, $\cal B$, and $\cal D$,
$${\cal A}\boxtimes ({\cal B}\boxtimes {\cal D}) \ne ({\cal A}\boxtimes {\cal B})\boxtimes {\cal D}$$
in general. Besides, this $\boxtimes$ product is different from other types of existing tensor products in the literature \cite{KB, QL}.

\begin{0}
For $\alpha \ge 2$, the $\alpha$th power of an $m$th order, $n$ dimensional tensor $\cal A$, denoted as ${\cal A}^\alpha=[a^{(\alpha)}_{i_1i_2\ldots i_m}]$, is determined recursively by 
\be
\label{ak}
{\cal A}^\alpha={\cal A}^{\alpha-1} \boxtimes {\cal A}
\ee
with ${\cal A}^1={\cal A}$. Moreover, by convention, ${\cal A}^0$ is the $m$th order, $n$ dimensional identity tensor ${\cal I}=[\delta_{i_1i_2\ldots i_m}]$, whose entries are given by $\delta_{i_1i_2\ldots i_m}=\delta_{i_1i_2}$ for any $i_1, i_2, \ldots, i_m \in \langle n \rangle$, where, and in the sequel, $\delta_{i_1i_2}$ denotes the Kronecker delta, i.e.,
$$\delta_{i_1i_2}=\left\{\begin{array}{cl}
1, & i_1=i_2;\\
0, & {\rm otherwise}.
\end{array}\right.$$
\end{0}

Note that ${\cal I}$ is merely the left identity tensor, meaning that ${\cal I}\boxtimes {\cal A}={\cal A}$ but usually ${\cal A}\boxtimes {\cal I} \ne {\cal A}$. In addition, such an identity tensor is different from another existing one in the literature \cite{KB, QL}.

The preceding $\boxtimes$ product, power, and identity tensor have recently been implemented in MATLAB as {\tt bprod}, {\tt bpow}, and {\tt eyet}, respectively; see \cite{Xu26c}.

The $\boxtimes$ product and power have played a significant part in investigating higher order Markov chains, see \cite{HX24a, HX26, Xu26a, Xu26c}. Such a role is fundamentally due to a practical probabilistic interpretation as follows: When $\cal A$ stands for the transition probability tensor of an $(m-1)$th order Markov chain with $n$ states, ${\cal A}^\alpha$ turns out to be the $\alpha$-step transition tensor, consisting of all the $\alpha$-step transition probabilities of the chain. In particular, ${\cal A}^\alpha$ has been used to introduce \cite{HWX, HX24a, HX26}:

\begin{0}
\label{erg}
Let ${\cal A}$ be the transition tensor of an $(m-1)$th order, $n$-state Markov chain. Denote ${\cal A}^\alpha=[a^{(\alpha)}_{i_1i_2\ldots i_m}]$. If for any $i_1, i_2, \ldots, i_m \in \langle n \rangle$, there exists $\alpha \ge 1$, which may depend on $i_1, i_2, \ldots, i_m$, so that $a^{(\alpha)}_{i_1i_2\ldots i_m}>0$, then this chain (and $\cal A$ as well) is called ergodic.
\end{0}

\begin{0}
\label{reg}
Let $\cal A$ be the transition tensor of an $(m-1)$th order, $n$-state Markov chain. Then, this chain (and $\cal A$ too) is called regular if there exists $\alpha \ge 1$ such that ${\cal A}^\alpha>0$.
\end{0}

Extending Definitions \ref{erg} and \ref{reg} to nonnegative $m$th order, $n$ dimensional tensors, we now arrive at: 

\begin{0}
\label{irr}
Given a nonnegative $m$th order, $n$ dimensional tensor ${\cal A}$, denote ${\cal A}^\alpha=[a^{(\alpha)}_{i_1i_2\ldots i_m}]$. If for any $i_1, i_2, \ldots, i_m \in \langle n \rangle$, there exists $\alpha \ge 1$, which may depend on $i_1, i_2, \ldots, i_m$, such that $a^{(\alpha)}_{i_1i_2\ldots i_m}>0$, then $\cal A$ is called s-irreducible.
\end{0}

\begin{0}
\label{pri}
Given a nonnegative $m$th order, $n$ dimensional tensor ${\cal A}$, if there exists $\alpha \ge 1$ such that ${\cal A}^\alpha>0$, then $\cal A$ is called s-primitive.
\end{0}

Like the relationship between irreducibility and primitivity, it is obvious that s-primitivity implies s-irreducibility, but the converse does not hold. In addition, for the special case of $m=2$, i.e., for nonnegative $n \times n$ matrices, Definitions \ref{irr} and \ref{pri} clearly coincide with their respective counterparts in Definitions \ref{irr1} and \ref{pri1}. As we shall demonstrate next, for the general $m \ge 3$ tensor case, Definitions \ref{irr} and \ref{pri} manifest different concepts as compared with Definitions \ref{irr3} and \ref{pri3}, respectively.

Before continuing, we observe that by (\ref{ak}), 
\be
\label{aka}
a^{(\alpha)}_{i_1i_2\ldots i_m}=\sum_{j_2, \ldots, j_\alpha \in \langle n \rangle}\underbrace{a_{i_1j_2\ldots j_m}a_{j_2\ldots j_{m+1}}\cdots a_{j_{\alpha-m+1}\ldots j_\alpha}a_{j_{\alpha-m+2}\ldots j_\alpha i_2}\cdots a_{j_\alpha i_2 \ldots i_m}}_\alpha
\ee
for $\alpha > m$. The case $2 \le \alpha \le m$ leads to a similar expression, in which the summation takes the form 
\be
\label{akam}
a^{(\alpha)}_{i_1i_2\ldots i_m}=\sum_{j_2, \ldots, j_\alpha \in \langle n \rangle}\underbrace{a_{i_1j_2\ldots j_\alpha i_2\ldots i_{m-\alpha+1}}a_{j_2\ldots j_\alpha i_2\ldots i_{m-\alpha+2}} \cdots a_{j_\alpha i_2\ldots i_m}}_\alpha.
\ee
The s-irreducibility of the nonnegative tensor $\cal A$, therefore, means that once $i_1, i_2, \ldots, i_m$ are given, one of the terms in (\ref{aka}) or (\ref{akam}) is positive for some $\alpha \ge 1$. In addition, this $\alpha$ may vary for different $i_1, i_2, \ldots, i_m$. On the other hand, the s-primitivity of the nonnegative tensor $\cal A$ means that one value of $\alpha \ge 1$ works for all $i_1, i_2, \ldots, i_m$ such that one of the terms in (\ref{aka}) or (\ref{akam}) is positive. As a result, nonnegative s-irreducible and s-primitive tensors are not uncommon.

\begin{6}
\label{irrpriex0}
Take a third order, $4$ dimensional tensor 
$${\cal A}(:,:,1)=\left[\begin{array}{cccc}
1 & 0 & 0 & 0\\
1 & 0 & 1 & 0\\
0 & 1 & 0 & 1\\
0 & 0 & 0 & 0
\end{array}\right], ~{\cal A}(:,:,2)=\left[\begin{array}{cccc}
0 & 0 & 1 & 1\\
0 & 1 & 0 & 0\\
0 & 1 & 0 & 0\\
1 & 0 & 1 & 0
\end{array}\right],$$
$${\cal A}(:,:,3)=\left[\begin{array}{cccc}
0 & 1 & 0 & 1\\
1 & 0 & 1 & 0\\
0 & 0 & 1 & 0\\
0 & 0 & 0 & 0
\end{array}\right], ~{\cal A}(:,:,4)=\left[\begin{array}{cccc}
0 & 0 & 0 & 0\\
1 & 1 & 1 & 0\\
0 & 0 & 0 & 1\\
0 & 0 & 0 & 1
\end{array}\right].$$
We observe ${\cal A}+{\cal A}^2+{\cal A}^3+{\cal A}^4 > 0$. This implies that for any $i_1, i_2, i_3 \in \langle 4 \rangle$, there exists some $1 \le \alpha \le 4$, which turns out to depend on $i_1, i_2, i_3$, such that $a^{(\alpha)}_{i_1i_2i_3}>0$. Thus, $\cal A$ is s-irreducible. This $\cal A$ will be revisited in Example \ref{irrpriex1} after new results are developed. We shall see that $\cal A$ here is not s-primitive, but it is primitive.
\end{6}

Incidentally, from now on, $a^{(\alpha)}_{i_1i_2\ldots i_m}$ always represents the $(i_1, i_2, \ldots, i_m)$th entry of some ${\cal A}^\alpha$.

\begin{6}
Let $\cal A$ be the following third order, $4$ dimensional tensor:
$${\cal A}(:,:,1)=\left[\begin{array}{cccc}
1 & 0 & 0 & 0\\
1 & 0 & 1 & 0\\
0 & 1 & 0 & 1\\
0 & 0 & 0 & 0
\end{array}\right], ~{\cal A}(:,:,2)=\left[\begin{array}{cccc}
0 & 0 & 1 & 1\\
0 & 1 & 0 & 0\\
1 & 1 & 0 & 0\\
1 & 0 & 1 & 0
\end{array}\right],$$
$${\cal A}(:,:,3)=\left[\begin{array}{cccc}
0 & 1 & 0 & 1\\
1 & 0 & 1 & 0\\
0 & 0 & 1 & 0\\
0 & 0 & 0 & 0
\end{array}\right], ~{\cal A}(:,:,4)=\left[\begin{array}{cccc}
0 & 0 & 0 & 0\\
1 & 1 & 1 & 0\\
0 & 0 & 0 & 1\\
0 & 0 & 0 & 1
\end{array}\right].$$
This $\cal A$ is obtained from the previous example by adding one positive entry. Even though nearly two thirds of the entries of $\cal A$ are zeros, we have ${\cal A}^{10} > 0$ and therefore $\cal A$ is s-primitive.
\end{6}

Regarding the relationship between irreducibility and s-irreducibility, we have:

\begin{2}
\label{s2irr}
For a nonnegative $m$th order, $n$ dimensional tensor $\cal A$, if $\cal A$ is s-irreducible, then it is irreducible.
\end{2}
\bp
See \cite[Theorem 3.2]{HWX}, whose proof is done on stochastic tensors, but can be easily extended to nonnegative $m$th order, $n$ dimensional tensors. 
\ep

The converse of Theorem \ref{s2irr} is not the case in general, as seen from the example below \cite{HWX}. In other words, Definitions \ref{irr} and \ref{irr3} are not equivalent for nonnegative tensors, although Definitions \ref{irr1} and \ref{irr2} are equivalent for nonnegative matrices.

\begin{6}
Consider a third order, $3$ dimensional tensor
$${\cal A}(:,:,1)=\left[\begin{array}{ccc}
0 & 0 & 0\\
1 & 0 & 0\\
0 & 1 & 1
\end{array}\right], ~{\cal A}(:,:,2)=\left[\begin{array}{ccc}
0 & 0 & 0\\
0 & 0 & 0\\
1 & 1 & 1
\end{array}\right], ~{\cal A}(:,:,3)=\left[\begin{array}{ccc}
0 & 0 & 1\\
0 & 0 & 0\\
1 & 1 & 0
\end{array}\right].$$
It is easy to confirm that $\cal A$ is irreducible, yet $\cal A$ is not s-irreducible because $a^{(\alpha)}_{2i_2i_3}=0$ for any $i_2, i_3 \in \langle 3 \rangle$ and any $\alpha \ge 2$.
\end{6}

Another immediate observation regarding s-primitivity is given below.

\begin{2}
Let $\cal A$ be a nonnegative $m$th order, $n$ dimensional tensor such that ${\cal A}^\alpha > 0$ for some $\alpha \ge 1$. Then ${\cal A}^\beta > 0$ for any $\beta > \alpha$.
\end{2}
\bp
It suffices to show that ${\cal A}^{\alpha+1} > 0$. Observe that, according to (\ref{bprod}), $c_{i_1i_2\ldots i_m}=0$ if $b_{ji_2\ldots i_m}=0$ for all $j \in \langle n \rangle$. Assume, to the contrary, that 
$$a^{(\alpha+1)}_{i_1i_2\ldots i_m}=\sum_{j \in \langle n \rangle}a^{(\alpha)}_{i_1ji_2\ldots i_{m-1}}a_{ji_2\ldots i_m}=0$$
for some $i_1, i_2, \ldots, i_m \in \langle n \rangle$. Then, we must have $a_{ji_2\ldots i_m}=0$ for all $j \in \langle n \rangle$. This, however, leads to $a^{(\alpha)}_{i_1i_2\ldots i_m}=0$ for any $\alpha \ge 1$, thus a contradiction. 
\ep

The relationship between primitivity and s-primitivity appears to be not that straightforward in a manner analogous to Theorem \ref{s2irr}. To begin with, nevertheless, we have the following special case, whose proof is trivial and is thus omitted.

\begin{2}
Any positive $m$th order, $n$ dimensional tensor is both s-primitivity and primitivity.
\end{2}

Next, let us quote a result in \cite{CPZ} as a technical lemma.

\begin{1}
\label{incr}
Given a nonnegative $m$th order, $n$ dimensional tensor $\cal A$, let $T_{\cal A}$ be defined as (\ref{map}). Then, $T_{\cal A}$ is nondecreasing, i.e., for any $n$-vectors $x$ and $y$ such that $y \ge x \ge 0$, $T_{\cal A}(y) \ge T_{\cal A}(x)$.
\end{1}

Moreover, we can state the lemma below, whose proof is straightforward and is therefore omitted.

\begin{1}
\label{add}
Given a nonnegative $m$th order, $n$ dimensional tensor $\cal A$, let $T_{\cal A}$ be defined as (\ref{map}). Then, $T_{\cal A}(x+y) \ge T_{\cal A}(x)+T_{\cal A}(y)$ for any nonnegative $n$-vectors $x$ and $y$.
\end{1}

To help address the question how s-primitivity is different from primitivity, we now develop two characterizations of the latter. These characterizations are useful in checking primitivity and constructing primitive tensors. The second one also highlights how primitivity may differ from s-primitivity. 

\begin{2}
\label{te}
Let $\cal A$ be a nonnegative $m$th order, $n$ dimensional tensor. Then, $\cal A$ is primitive if and only if there exists $\alpha \ge 1$ such that $T_{\cal A}^{\,\circ \alpha}(e_i) > 0$ for all $i \in \langle n \rangle$, where $T_{\cal A}$ is given by (\ref{map}).
\end{2}
\bp
The necessity is obvious by Definition \ref{pri3}. To prove the sufficiency, we take any $n$-vector $0 \ne x \ge 0$. Observe that $x \ge ce_i$ for some scalar $c > 0$ and $i \in \langle n \rangle$. Using Lemma \ref{incr}, we have 
$$T_{\cal A}(x) \ge c^{m-1}T_{\cal A}(e_i),$$
$$T_{\cal A}^{\,\circ 2}(x) \ge c^{2(m-1)}T_{\cal A}^{\,\circ 2}(e_i),$$
$$\vdots$$
$$T_{\cal A}^{\,\circ \alpha}(x) \ge c^{\alpha (m-1)}T_{\cal A}^{\,\circ \alpha}(e_i) > 0.$$
This completes the proof.
\ep

\begin{2}
\label{jiii}
Let $A$ be a nonnegative $m$th order, $n$ dimensional tensor and let $T_{\cal A}$ be given as in (\ref{map}). Then, $\cal A$ is primitive if there exists $\alpha \ge 1$ such that for any $i, j_\alpha \in \langle n \rangle$, 
\be
\label{chain}
\underbrace{a_{j_1i\ldots i}a_{j_2j_1\ldots j_1}a_{j_3j_2\ldots j_2}\cdots a_{j_\alpha j_{\alpha-1}\ldots j_{\alpha-1}}}_\alpha > 0
\ee
for some $j_1, j_2, \ldots, j_{\alpha-1} \in \langle n \rangle$, which may depend on $i$ and $j_\alpha$. In particular, when $\alpha=1$, (\ref{chain}) is interpreted as $a_{j_1i\ldots i} > 0$ for any $i, j_1 \in \langle n \rangle$.
\end{2}
\bp
In view of Theorem \ref{te}, it suffices for us to consider $T_{\cal A}^{\,\circ \alpha}(e_i) > 0$ for some $\alpha \ge 1$ and any $i \in \langle n \rangle$.

When $\alpha=1$, we see
$$T_{\cal A}(e_i)=\sum_{j_1 \in \langle n \rangle} a_{j_1i\ldots i}e_{j_1}.$$
Hence, $T_{\cal A}(e_i) > 0$ for all $i \in \langle n \rangle$ if and only if $a_{j_1i\ldots i}>0$ for any $i, j_1 \in \langle n \rangle$.

When $\alpha=2$, it follows from the above and Lemma \ref{add} that 
$$T_{\cal A}^{\,\circ 2}(e_i) \ge \sum_{j_1 \in \langle n \rangle} a_{j_1i\ldots i}^{m-1}T_{\cal A}(e_{j_1})=\sum_{j_2 \in \langle n \rangle}\sum_{j_1 \in \langle n \rangle} a_{j_1i\ldots i}^{m-1}a_{j_2j_1\ldots j_1}e_{j_2},$$
hence $T_{\cal A}^{\,\circ 2}(e_i) > 0$ for all $i \in \langle n \rangle$ if for any $i, j_2 \in \langle n \rangle$, $a_{j_1i\ldots i}^{m-1}a_{j_2j_1\ldots j_1} > 0$ for some $j_1 \in \langle n \rangle$; or equivalently, 
$a_{j_1i\ldots i}a_{j_2j_1\ldots j_1} > 0$ for some $j_1 \in \langle n \rangle$ thanks to the nonnegativity of $\cal A$. Clearly, $j_1$ here may depend on $i$ and $j_2$.

In general, for any $\alpha \ge 2$, we have 
$$T_{\cal A}^{\,\circ \alpha}(e_i) \ge \sum_{j_\alpha \in \langle n \rangle}\sum_{j_{\alpha-1} \in \langle n \rangle}{\underbrace{\left(\ldots \left(\sum_{j_2 \in \langle n \rangle}\left(\sum_{j_1 \in \langle n \rangle} a_{j_1i\ldots i}^{m-1}a_{j_2j_1\ldots j_1}\right)^{m-1}\!\!\!\!\!a_{j_3j_2\ldots j_2}\right)^{m-1}\!\!\!\!\!\ldots \right)}_{\alpha-2}}^{m-1}\!\!\!\!\!a_{j_\alpha j_{\alpha-1}\ldots j_{\alpha-1}}e_{j_\alpha}.$$
implying that $T_{\cal A}^{\,\circ \alpha}(e_i) > 0$ for all $i \in \langle n \rangle$ if for any $i, j_\alpha \in \langle n \rangle$,
$$\sum_{j_{\alpha-1} \in \langle n \rangle}{\underbrace{\left(\ldots \left(\sum_{j_2 \in \langle n \rangle}\left(\sum_{j_1 \in \langle n \rangle} a_{j_1i\ldots i}^{m-1}a_{j_2j_1\ldots j_1}\right)^{m-1}a_{j_3j_2\ldots j_2}\right)^{m-1}\ldots \right)}_{\alpha-2}}^{m-1}a_{j_\alpha j_{\alpha-1}\ldots j_{\alpha-1}} > 0.$$
Again, due to the nonnegativity of $\cal A$, this latter inequality is guaranteed if for any $i, j_\alpha \in \langle n \rangle$, (\ref{chain}) holds for some $j_1, j_2, \ldots, j_{\alpha-1} \in \langle n \rangle$, where these $j_1, j_2, \ldots, j_{\alpha-1}$ may depend on $i$ and $j_\alpha$.
\ep

Although it gives only a sufficient condition for primitivity, Theorem \ref{jiii} shows a main distinction between primitivity and s-primitivity. Specifically, while s-primitivity hinges on some $\boxtimes$ power involving, in general, the entire nonnegative tensor $\cal A$, primitivity may be accomplished via the entries of $n$ mode-$1$ columns (fibers) ${\cal A}(:,i,i,\ldots,i)$, $i \in \langle n \rangle$. It is from this perspective that s-primitivity tends to be a stronger condition than primitivity.

We are now in a position to give two examples. From these examples, in particular, we see that primitivity does not necessarily imply s-irreducibility and that s-primitivity does not necessarily imply primitivity --- also refer to Example \ref{irrpriex1}, where the tensor is both primitive and s-irreducible.

\begin{6}
\label{priex1}
Consider a fourth order, $3$ dimensional tensor, which is given by 
$${\cal A}(:,:,1,1)=\left[\begin{array}{ccc}
1 & 0 & 0\\
0 & 0 & 0\\
1 & 0 & 0
\end{array}\right], ~{\cal A}(:,:,2,2)=\left[\begin{array}{ccc}
0 & 1 & 0\\
0 & 1 & 0\\
0 & 0 & 0
\end{array}\right], ~{\cal A}(:,:,3,3)=\left[\begin{array}{ccc}
0 & 0 & 0\\
0 & 0 & 1\\
0 & 0 & 1
\end{array}\right],$$
with the other frontal slices being $3 \times 3$ zero matrices. Notice that $\cal A$ has positive entries only in three mode-$1$ columns ${\cal A}(:,i,i,i)$, $i \in \langle 3 \rangle$. It is straightforward to verify that $\cal A$ is irreducible. Meanwhile, we have $T_{\cal A}^{\,\circ 2}(e_i)>0$ for any $i \in \langle 3 \rangle$, which implies the primitivity of $\cal A$. This $\cal A$, however, is neither s-irreducible nor s-primitive since ${\cal A}^\alpha={\cal A}$ for all $\alpha \ge 1$.
\end{6}

\begin{6}
\label{priex2}
Take a fourth order, $3$ dimensional tensor $\cal A$ such that
$${\cal A}(:,:,1,1)=\left[\begin{array}{ccc}
0 & 1 & 1\\
0 & 1 & 1\\
1 & 1 & 1
\end{array}\right], ~{\cal A}(:,:,2,2)=\left[\begin{array}{ccc}
1 & 1 & 1\\
1 & 0 & 1\\
1 & 0 & 1
\end{array}\right], ~{\cal A}(:,:,3,3)=\left[\begin{array}{ccc}
1 & 1 & 0\\
1 & 1 & 1\\
1 & 1 & 0
\end{array}\right],$$
and all the remaining frontal slices are $3 \times 3$ matrices of all ones. It is easy to check that $\cal A$ is s-primitive since ${\cal A}^2 > 0$. On the other hand, we observe 
$$T_{\cal A}(e_1)=e_3, ~T_{\cal A}(e_2)=e_1, ~T_{\cal A}(e_3)=e_2.$$
Because of this, for each $i \in \langle 3 \rangle$, $T_{\cal A}^{\,\circ \alpha}(e_i)$ rotates through $e_1, e_2, e_3$ as $\alpha$ increases. Hence, $\cal A$ is not primitive. Meanwhile, since $\cal A$ is s-primitive, it is s-irreducible and thus also irreducible.
\end{6}

Examples \ref{priex1} and \ref{priex2} also demonstrate the key role the mode-$1$ columns ${\cal A}(:,i,i,\ldots,i)$ may play in determining primitivity. In Example \ref{priex1}, for instance, if we check $T_{\cal A}^{\,\circ 2}(e_1)$, then $a_{j_1 111}a_{1j_1j_1j_1} > 0$ at $j_1=1$, $a_{j_1 111}a_{2j_1j_1j_1} > 0$ at $j_1=3$, and $a_{j_1 111}a_{3j_1j_1j_1} > 0$ at $j_1=3$. Thus, condition (\ref{chain}) is satisfied when $\alpha=2$ and $i=1$. Similarly, we can easily verify that condition (\ref{chain}) is also satisfied when $\alpha=2$ and $i=2, 3$. For primitivity alone, the remaining entries of $\cal A$ are no longer relevant once we have condition (\ref{chain}) at hand. In Example \ref{priex2}, on the other hand, condition (\ref{chain}) cannot be met because the only product of positive entries from those mode-$1$ columns must involve a contiguous segment of 
$$\cdots a_{3111}a_{2333}a_{1222}a_{3111}a_{2333}a_{1222}\cdots.$$
Given $i=1$ and any fixed $\alpha \ge 1$, for example, it is just impossible to realize (\ref{chain}) for all $j_\alpha \in \langle 3 \rangle$.

In the remainder of this work, we shall further explore the relationship between s-irreducibility and s-primitivity. For this purpose, we need several preparatory results.

The set of multi-indices of length $\alpha \ge 2$, denoted by $\langle n \rangle^\alpha$, consists of all the indices in the form $i_1i_2\ldots i_\alpha$ such that $i_j \in \langle n \rangle$ for $j=1, 2, \ldots, \alpha$. For convenience, we shall always order these multi-indices using linear indexing \cite{MSL}. Denote $N=n^\alpha$. Then, the entries of an $N \times N$ matrix $Q$ can also be written in multi-index form as 
$$q_{i_1i_2\ldots i_\alpha,j_1j_2\ldots j_\alpha}, ~i_1i_2\ldots i_\alpha, ~j_1j_2\ldots j_\alpha \in \langle n \rangle^\alpha.$$
The next two definitions generalize the notions of the transition matrix of a reduced first order Markov chain and matricization of the transition tensor of a higher order Markov chain \cite{HX24a, HX26, KB, Xu26a, Xu26b}.

\begin{0}
\label{rcmat}
For an $m$th order, $n$ dimensional tensor $\cal A$, its reduced matrix $Q_{\cal A}$ is of size $N \times N$, where $N=n^{m-1}$, and its entries are such that for any $i_1i_1\ldots, i_{m-1}, \,j_2j_3\ldots j_m \in \langle n \rangle^{m-1}$, 
\be
\label{q}
q_{i_1i_2\ldots i_{m-1},j_2j_3\ldots j_m}=\left\{\begin{array}{cl}
a_{i_1i_2\ldots i_{m-1}j_m}, & i_\ell=j_\ell, ~\ell=2, 3, \ldots, m-1;\\
0, & {\rm otherwise.}
\end{array}\right.
\ee
\end{0}

\begin{0}
\label{matri}
For an $m$th order, $n$ dimensional tensor $\cal A$, its mode-$1$ matricization $A$ is an $n \times N$ matrix, where $N=n^{m-1}$, whose entries are given, in multi-index form, by 
$$a_{i_1,i_2i_3\ldots i_m}, ~i_1 \in \langle n \rangle, ~i_2\ldots i_m \in \langle n \rangle^{m-1}.$$
In other words, this matrix $A$ is formed by arranging the frontal slices of $\cal A$ side by side in the linear indexing order of $i_3, \ldots, i_m$.
\end{0}

The mode-$1$ matricization of ${\cal A}^\alpha$ will be denoted accordingly as $A^{(\alpha)}$. In particular, the entries of $A^{(0)}$ are simply $a^{(0)}_{i_1,i_2\ldots i_m}=\delta_{i_1i_2}$. We mention that the validity of this expression of $a^{(0)}_{i_1,i_2\ldots i_m}$ hinges on linear indexing order.

One way of constructing $Q_{\cal A}$ is given by \cite{HX26} 
$$Q_{\cal A}=G\ast A,$$
where $G=[\underbrace{I_{n^{m-2}} \ \ \ldots \ \ I_{n^{m-2}}}_{n}]$ and $\ast$ is the Khatri-Rao product \cite{SBG}, which has been implemented along with mode-$1$ matricization as MATLAB functions {\tt rcmat} and {\tt t2mat}, respectively, in \cite{Xu26c}.

As a by-product, using the mode-$1$ matricization, the map in (\ref{map}) can be written as 
$$T_{\cal A}(x)=A(\underbrace{x\otimes_K \cdots \otimes_K x}_{m-1}),$$
where $\otimes_K$ denotes the Kronecker product and the entries of $\underbrace{x\otimes_K \cdots \otimes_K x}_{m-1}$, i.e., $x_{i_2}\ldots x_{i_m}$, may be arranged via the linear indexing order of $i_m\ldots i_2$.

In passing, we provide an example below to illustrate Definitions \ref{rcmat} and \ref{matri}.

\begin{6}
Let ${\cal A}$ be a third order, $3$ dimensional tensor. Then 
$$Q_{\cal A}=\left[\ba{ccccccccc}
a_{111} & 0 & 0 & a_{112} & 0 & 0 & a_{113} & 0 & 0\\
a_{211} & 0 & 0 & a_{212} & 0 & 0 & a_{213} & 0 & 0\\
a_{311} & 0 & 0 & a_{312} & 0 & 0 & a_{313} & 0 & 0\\
0 & a_{121} & 0 & 0 & a_{122} & 0 & 0 & a_{123} & 0\\
0 & a_{221} & 0 & 0 & a_{222} & 0 & 0 & a_{223} & 0\\
0 & a_{321} & 0 & 0 & a_{322} & 0 & 0 & a_{323} & 0\\
0 & 0 & a_{131} & 0 & 0 & a_{132} & 0 & 0 & a_{133}\\
0 & 0 & a_{231} & 0 & 0 & a_{232} & 0 & 0 & a_{233}\\
0 & 0 & a_{331} & 0 & 0 & a_{332} & 0 & 0 & a_{333}
\ea\right]$$
and
$$A=\left[\begin{array}{ccccccccc}
a_{111} & a_{121} & a_{131} & a_{112} & a_{122} & a_{132} & a_{113} & a_{123} & a_{133}\\
a_{211} & a_{221} & a_{231} & a_{212} & a_{222} & a_{232} & a_{213} & a_{223} & a_{233}\\
a_{311} & a_{321} & a_{331} & a_{312} & a_{322} & a_{332} & a_{313} & a_{323} & a_{333}
\end{array}\right].$$
\end{6}

For stochastic transition tensors of any order $m$, a probabilistic proof of the result below is presented in \cite{HX26}. The $m=3$ case can also be seen, with no proof, in \cite{Vla85}. We present here an algebraic proof.

\begin{2}
\label{komo}
Let $\cal A$ be an $m$th order, $n$ dimensional tensor. Then, for any $\alpha, \beta \ge 0$, 
\be
\label{aq}
A^{(\alpha+\beta)}=A^{(\alpha)}Q_{\cal A}^\beta.
\ee
\end{2}
\bp
Clearly, (\ref{aq}) is true whenever $\beta=0$. For $\beta=1$, we use (\ref{q}) and (\ref{bprod}) to obtain that for any $i_1, i_2, \ldots, i_m \in \langle n \rangle$,
\be
\begin{split}
(A^{(\alpha)}Q_{\cal A})_{i_1,i_2i_3\ldots i_m} &=\sum_{j_1j_2\ldots j_{m-1} \in \langle n \rangle^{m-1}} a^{(\alpha)}_{i_1,j_1j_2\ldots j_{m-1}}q_{j_1j_2\ldots j_{m-1},i_2i_3\ldots i_m}\\ \nonumber
&=\sum_{j_1 \in \langle n \rangle}a^{(\alpha)}_{i_1j_1i_2\ldots i_{m-1}}a_{j_1i_2\ldots i_m}\\
&=a^{(\alpha+1)}_{i_1,i_2\ldots i_m},
\end{split}
\ee
i.e., $A^{(\alpha)}Q_{\cal A}=A^{(\alpha+1)}$ for any $\alpha \ge 0$. Consequently, we see that for $\beta=2$, 
$$A^{(\alpha)}Q_{\cal A}^2=(A^{(\alpha)}Q_{\cal A})Q_{\cal A}=A^{(\alpha+1)}Q_{\cal A}=A^{(\alpha+2)}.$$
Proceeding in this fashion, the conclusion follows.
\ep

An immediate consequence of Theorem \ref{komo} is the following:

\begin{2}
\label{q2a}
For a nonnegative $m$th order, $n$ dimensional tensor $\cal A$, if its reduced matrix $Q_{\cal A}$ is primitive, then $\cal A$ is s-primitive.
\end{2}
\bp
See \cite[Theorem 3.1]{HX26}, whose proof is formulated on stochastic tensors. The same argument, however, also works on nonnegative $m$th order, $n$ dimensional tensors.
\ep

The reverse of Theorem \ref{q2a}, in general, is not true. To verify this, let us consider: 

\begin{6}
Choose a third order, $3$ dimensional tensor 
$${\cal A}(:,:,1)=\left[\begin{array}{ccc}
0 & 0 & 1\\
1 & 1 & 0\\
0 & 0 & 1
\end{array}\right], ~{\cal A}(:,:,2)=\left[\begin{array}{ccc}
1 & 1 & 0\\
0 & 1 & 0\\
1 & 0 & 1
\end{array}\right], ~{\cal A}(:,:,3)=\left[\begin{array}{ccc}
1 & 0 & 0\\
0 & 1 & 1\\
1 & 0 & 1
\end{array}\right].$$
Since ${\cal A}^4 > 0$, $\cal A$ is s-primitive. Meanwhile, we have 
$$Q_{\cal A}=\left[\begin{array}{ccccccccc}
0  &  0  &  0  &  1  &  0  &  0  &  1  &  0  &  0\\
1  &  0  &  0  &  0  &  0  &  0  &  0  &  0  &  0\\
0  &  0  &  0  &  1  &  0  &  0  &  1  &  0  &  0\\
0  &  0  &  0  &  0  &  1  &  0  &  0  &  0  &  0\\
0  &  1  &  0  &  0  &  1  &  0  &  0  &  1  &  0\\
0  &  0  &  0  &  0  &  0  &  0  &  0  &  0  &  0\\
0  &  0  &  1  &  0  &  0  &  0  &  0  &  0  &  0\\
0  &  0  &  0  &  0  &  0  &  0  &  0  &  0  &  1\\
0  &  0  &  1  &  0  &  0  &  1  &  0  &  0  &  1
\end{array}\right].$$
Clearly, $Q_{\cal A}$ is not even irreducible because it has a row of zeros. This shows that $Q_{\cal A}$ is not primitive.
\end{6}

For the special case where ${\cal A} > 0$, however, the s-primitivity of $\cal A$ does translate to the primitivity of $Q_{\cal A}$, which is stated below.

\begin{2}
If $\cal A$ is a positive $m$th order, $n$ dimensional tensor, then its reduced matrix $Q_{\cal A}$ must be primitive.
\end{2}
\bp
The proof can be done in the same fashion as that for \cite[Theorem 3.2]{HX26}. The proof there is carried out on stochastic tensors, but its argument can be easily adapted for general positive $m$th order, $n$ dimensional tensors.
\ep

Moving on, we shall also treat $\langle n \rangle$ as the set of vertices, nodes, or states in the rest of this work. Let us cite a notion of accessibility \cite{HX24b, Vla85} first.

\begin{0}
\label{acc}
Let ${\cal A}$ be a nonnegative $m$th order, $n$ dimensional tensor and let $i, j \in \langle n \rangle$. Then, $j$ is said to be accessible from $i$ in $\cal A$, written as $i \rightarrow j$, if for any $i_3, \ldots, i_m \in \langle n \rangle$, there exists $\alpha \ge 0$, which may depend on $i_3, \ldots, i_m$, such that $a^{(\alpha)}_{jii_3\ldots i_m}>0$.
\end{0}

Similar to the discussion following Definitions \ref{irr} and \ref{pri}, for $i \rightarrow j$, one of the terms in (\ref{aka}) or (\ref{akam}) --- with $i_1$ and $i_2$ there being replaced by $i$ and $j$ --- is positive for some $\alpha \ge 1$ and, additionally, this value of $\alpha$ may vary depending on $i_3, \ldots, i_m$. In particular, if $\cal A$ is s-irreducible, then $i \rightarrow j$ for any $i, j \in \langle n \rangle$. Such an accessibility condition, therefore, is not particularly demanding. 

Furthermore, we can specify $i \leftrightarrow j$ as both $i \rightarrow j$ and $j \rightarrow i$. Such a $\leftrightarrow$ relationship turns out to be an equivalence relationship \cite{HX24b, Vla85}. The reason that $\alpha$ is allowed to be zero in Definition \ref{acc} is essentially for the reflexivity of this $\leftrightarrow$ equivalence relationship \cite{HX24b}.

For a nonnegative $m$th order, $n$ dimensional tensor $\cal A$, we set 
\be
\label{sa}
S_{\cal A}(i,j)=\{\alpha \ge 1 : a^{(\alpha)}_{jii_3\ldots i_m}>0 {\rm ~for ~any} ~i_3, \ldots, i_m \in \langle n \rangle\}.
\ee
Notice that $S_{\cal A}(i,j) \ne \emptyset$ is a condition stronger than $i \rightarrow j$. By convention, $S_{\cal A}(i)=S_{\cal A}(i,i)$.

We are now ready to broaden a number of results regarding third order stochastic transition tensors in \cite{Vla84, Vla85} to nonnegative $m$th order, $n$ dimensional tensors.

\begin{1}
\label{sij}
Given a nonnegative $m$th order, $n$ dimensional tensor $\cal A$, let $i, j, k \in \langle n \rangle$. If $\alpha \in S_{\cal A}(i,j)$ and $\beta \in S_{\cal A}(j,k)$, then $\alpha+\beta \in S_{\cal A}(i,k)$.
\end{1}
\bp
By (\ref{aq}) in Theorem \ref{komo}, we know  
\be
\begin{split}
a^{(\alpha)}_{j,ii_3\ldots i_m} &=\sum_{j_1j_2\ldots j_{m-1} \in \langle n \rangle^{m-1}}a^{(0)}_{j,j_1j_2\ldots j_{m-1}}q^{(\alpha)}_{j_1j_2\ldots j_{m-1},ii_3\ldots i_m}\\ \nonumber
&=\sum_{j_1j_2\ldots j_{m-1} \in \langle n \rangle^{m-1}}\delta_{jj_1}q^{(\alpha)}_{j_1j_2\ldots j_{m-1},ii_3\ldots i_m}\\
&=\sum_{j_2\ldots j_{m-1} \in \langle n \rangle^{m-2}}q^{(\alpha)}_{jj_2\ldots j_{m-1},ii_3\ldots i_m} > 0.
\end{split}
\ee
Note that the latter inequality holds since $\alpha \in S_{\cal A}(i,j)$.

Next, we employ (\ref{aq}) again to obtain 
\be
\begin{split}
a^{(\alpha+\beta)}_{k,ii_3\ldots i_m} & =\sum_{j_1j_2\ldots j_{m-1} \in \langle n \rangle^{m-1}}a^{(\beta)}_{k,j_1j_2\ldots j_{m-1}}q^{(\alpha)}_{j_1j_2\ldots j_{m-1},ii_3\ldots i_m}\\ \nonumber
& \ge \sum_{j_2\ldots j_{m-1} \in \langle n \rangle^{m-2}}a^{(\beta)}_{k,jj_2\ldots j_{m-1}}q^{(\alpha)}_{jj_2\ldots j_{m-1},ii_3\ldots i_m}\\
& \ge \min_{j_2\ldots j_{m-1} \in \langle n \rangle^{m-2}}a^{(\beta)}_{k,jj_2\ldots j_{m-1}}\sum_{j_2\ldots j_{m-1} \in \langle n \rangle^{m-2}}q^{(\alpha)}_{jj_2\ldots j_{m-1},ii_3\ldots i_m} > 0
\end{split}
\ee
for any $i_3, \ldots, i_m \in \langle n \rangle$. In the last inequality, the fact $\beta \in S_{\cal A}(j,k)$ has also been utilized.
\ep

\begin{0}
\label{ape}
For a nonnegative $m$th order, $n$ dimensional tensor $\cal A$, $i \in \langle n \rangle$ is said to be aperiodic in $\cal A$ if  $gcd(S_{\cal A}(i))=1$, where, and in the sequel, gcd stands for the ``greatest common divisor''.
\end{0}

Incidentally, $d=\gcd(S_{\cal A}(i))$ is also called the period of $i$. If $S_{\cal A}(i)=\emptyset$, $d$ is regarded as zero. 

The following number theoretic conclusion can be found, for example, in \cite{BR}. This result can be formulated in terms of numerical semigroups \cite{RG} as well.

\begin{1}
\label{comb}
Let $\gamma_1, \gamma_2, \ldots, \gamma_k$ be positive integers such that $$\gcd \{\gamma_1, \gamma_2, \ldots, \gamma_k\}=1.$$ Then, there exists a positive integer $\beta$ such that for any positive integer $\alpha \ge \beta$, 
$$\alpha =\sum_{i=1}^k c_i\gamma_i$$
for some nonnegative integers $c_1, c_2, \ldots, c_k$.
\end{1}

The next result gives a full characterization of the relationship between s-irreducibility and s-primitivity. As mentioned earlier, a similar conclusion can be seen in \cite{Vla84, Vla85} for third order stochastic transition tensors. In these works, however, the role condition (iii) below plays has not been spelled out in full detail. Our proof here helps address this gap.

\begin{2}
\label{irrpri}
Consider a nonnegative $m$th order, $n$ dimensional tensor $\cal A$. Then, $\cal A$ is s-primitive if and only if: 
\begin{itemize}
\item[(i)]{$\cal A$ is s-irreducible,}
\item[(ii)]{$S_{\cal A}(i,j) \ne \emptyset$ whenever $i, j \in \langle n \rangle$ and $i \rightarrow j$ in $\cal A$, and}
\item[(iii)]{all $i \in \langle n \rangle$ are aperiodic in $\cal A$.}
\end{itemize}
\end{2}
\bp
The necessity part is straightforward. To see, for example, part (iii), we just resort to the fact that there exists $\beta \ge 1$ such that $a^{(\alpha)}_{iii_3\ldots i_m}>0$ for any $i, i_3, \ldots, i_m \in \langle n \rangle$ and $\alpha \ge \beta$. Thus, $\gcd(S_{\cal A}(i))=1$.

To prove the sufficiency part, we proceed as follows.

Given any $j, i, i_3, \ldots, i_m \in \langle n \rangle$, we need to show that there exists $\alpha \ge 1$, which is independent of $j, i, i_3, \ldots, i_m$, such that $a^{(\alpha)}_{jii_3\ldots i_m}>0$. Let us fix such arbitrary $j, i, i_3, \ldots, i_m$ first.

From parts (i) and (ii), we know $S_{\cal A}(i,j) \ne \emptyset$. Pick a fixed $\beta_{ij} \in S_{\cal A}(i,j)$. Similarly, we have $S_{\cal A}(j) \ne \emptyset$ and $\gcd(S_{\cal A}(j))=1$ by part (iii).

We claim now that there exists a finite subset $\tilde S_{\cal A}(j) \subseteq S_{\cal A}(j)$ such that $\gcd (\tilde S_{\cal A}(j))=1$. There is nothing to prove when $S_{\cal A}(j)$ is finite. Otherwise, we assume $S_{\cal A}(j)=\{\gamma_1, \gamma_2, \ldots\}$. Let
$$g_k=\gcd \{\gamma_1, \gamma_2, \ldots, \gamma_k\}, ~k=1, 2, \ldots.$$
It is clear that $\{g_k\}$ is a decreasing sequence and is also bounded below by $1$. Consequently, there are some positive integers $g, K \ge 1$ such that $g_k =g$ for any $k \ge K$. This, however, implies $g=\gcd(S_{\cal A}(j))=1$. Hence, it follows that we may choose $\tilde S_{\cal A}(j)=\{\gamma_1, \gamma_2, \ldots, \gamma_K\}$.

By Lemmas \ref{sij} and \ref{comb}, the above claim shows that there exists $\beta_j \ge 1$ such that $\alpha \in S_{\cal A}(j)$ for any $\alpha \ge \beta_j$. Meanwhile, from Lemma \ref{sij}, $\beta_{ij}+\beta_j \in S_{\cal A}(i,j)$. Next, we set 
$$\alpha=\max_{i,j \in \langle n \rangle}(\beta_{ij}+\beta_j).$$
Then, $\alpha-\beta_{ij} \ge \beta_j$ for any $i, j \in \langle n \rangle$, thus $\alpha-\beta_{ij} \in S_{\cal A}(j)$ for any $i, j \in \langle n \rangle$. Using, again, Lemma \ref{sij}, we conclude $\alpha=\beta_{ij}+(\alpha-\beta_{ij}) \in S_{\cal A}(i,j)$, which yields $a^{(\alpha)}_{jii_3\ldots i_m}>0$ for any $j, i, i_3, \ldots, i_m \in \langle n \rangle$. Hence, $\cal A$ is s-primitive.
\ep

In particular, for the special case of matrices, i.e., $m=2$, we notice that condition (ii) is automatically satisfied since (\ref{sa}) becomes 
$$S_{\cal A}(i,j)=\{\alpha \ge 1 : a^{(\alpha)}_{ji} > 0\},$$
which is nonempty whenever $i \rightarrow j$. Theorem \ref{irrpri}, therefore, reduces to the well-known classical result, see \cite{HJ} for example, as follows.

\begin{3}
A nonnegative $n \times n$ matrix is primitive if and only if it is irreducible and all $i \in \langle n \rangle$ are aperiodic in $A$.
\end{3}

To help illustrate Theorem \ref{irrpri}, we now revisit Example \ref{irrpriex0}.

\begin{6}
\label{irrpriex1}
Choose the same third order, $4$ dimensional s-irreducible tensor $\cal A$ as in Example \ref{irrpriex0}. Let us use Theorem \ref{irrpri} to decide whether $\cal A$ is s-primitive. For each $i \in \langle 4 \rangle$, we get $S_{\cal A}(i)=\{4, 8, 12, \ldots\}$, i.e., none of the vertices is aperiodic, which confirms that $\cal A$ is not s-primitive. Furthermore, in the same manner as in Examples \ref{priex1} and \ref{priex2}, it can be shown that $\cal A$ is primitive. Note that $\cal A$ here is both s-irreducible and primitive.
\end{6}

For any $m$th order, $n$ dimensional tensor ${\cal A}=[a_{i_1i_2\ldots i_m}]$, the associated diagonal tensor ${\cal A}_d=[a^{[d]}_{i_1i_2\ldots i_m}]$ is defined by 
$$a^{[d]}_{i_1i_2\ldots i_m}=a_{i_1i_2\ldots i_m}\delta_{i_1i_2},$$
where $i_1, i_2, \ldots, i_m \in \langle n \rangle$. Besides, $a_{iii_3\ldots i_m}$, $i, i_3, \ldots, i_m \in \langle n \rangle$, are called the diagonal entries of $\cal A$.

Continuing, let us state our last technical lemma. Its validity is clear and hence its proof is omitted.
\begin{1}
\label{comp}
Let $\cal A$ and $\cal B$ be both nonnegative $m$th order, $n$ dimensional tensors such that ${\cal A} \ge {\cal B}$. If $\cal B$ is s-irreducible, then so is $\cal A$. Similarly, if $\cal B$ is s-primitive, then so is $\cal A$.
\end{1}

Finally, we establish two more results concerning s-irreducibility and s-primitivity, both generalizing the corresponding classical results on nonnegative matrices \cite{BP, HJ}.

\begin{2}
\label{irr2pri}
For a nonnegative $m$th order, $n$ dimensional tensor $\cal A$, if it is s-irreducible, then $c_1{\cal I}+c_2{\cal A}$ is s-primitive, where $c_1, c_2$ are positive scalars and $\cal I$ is the identity tensor of the same size as $\cal A$.
\end{2}
\bp
We shall deal with the case when $c_1=c_2=1$ since the general case of $c_1, c_2 >0$ can be done in a similar fashion.

First, let us show that for all $\beta \ge 1$, $\ds ({\cal I}+{\cal A})^\beta \ge \sum_{\alpha=0}^\beta {\cal A}^\alpha$. Obviously, this holds trivially when $\beta=1$. For brevity, we denote 
$${\cal B}={\cal I}+{\cal A}.$$
When $\beta=2$, we have 
\be
\begin{split}
b^{(2)}_{i_1i_2\ldots i_m} &= \sum_{j \in \langle n \rangle}(\delta_{i_1j}+a_{i_1ji_2\ldots i_{m-1}})(\delta_{ji_2}+a_{ji_2\ldots i_m})\\ \nonumber
 &=\delta_{i_1i_2}+a_{i_1i_2\ldots i_m}+a^{(2)}_{i_1i_2\ldots i_m}+a_{i_1i_2i_2\ldots i_{m-1}}\\
 & \ge \delta_{i_1i_2}+a_{i_1i_2\ldots i_m}+a^{(2)}_{i_1i_2\ldots i_m},
\end{split}
\ee
i.e., $\ds ({\cal I}+{\cal A})^2 \ge \sum_{\alpha=0}^2 {\cal A}^\alpha$. Similarly, when $\beta=3$, we arrive at 
\be
\begin{split}
b^{(3)}_{i_1i_2\ldots i_m} &\ge \sum_{j \in \langle n \rangle}(\delta_{i_1j}+a_{i_1ji_2\ldots i_{m-1}}+a^{(2)}_{i_1ji_2\ldots i_{m-1}})(\delta_{ji_2}+a_{ji_2\ldots i_m})\\ \nonumber
 &=\delta_{i_1i_2}+a_{i_1i_2\ldots i_m}+a^{(2)}_{i_1i_2\ldots i_m}+a^{(3)}_{i_1i_2\ldots i_m}+a_{i_1i_2i_2\ldots i_{m-1}}+a^{(2)}_{i_1i_2i_2\ldots i_{m-1}}\\
 & \ge \delta_{i_1i_2}+a_{i_1i_2\ldots i_m}+a^{(2)}_{i_1i_2\ldots i_m}+a^{(3)}_{i_1i_2\ldots i_m},
\end{split}
\ee
i.e., $\ds ({\cal I}+{\cal A})^3 \ge \sum_{\alpha=0}^3 {\cal A}^\alpha$. The general case for any $\beta \ge 1$ follows by repeating the preceding argument.

Next, by the s-irreducibility of $\cal A$, for any $i_1, i_2, \ldots, i_m \in \langle n \rangle$, there exists $\alpha \ge 1$, which may depend on $i_1, i_2, \ldots, i_m$, such that $a^{(\alpha)}_{i_1i_2\ldots i_m} > 0$. Set 
$$\beta=\max_{i_1, i_2, \ldots, i_m \in \langle n \rangle}\min\{\alpha \ge 1 : a^{(\alpha)}_{i_1i_2\ldots i_m} > 0\}.$$
This implies that, for any $i_1, i_2, \ldots, i_m \in \langle n \rangle$, $a^{(\alpha)}_{i_1i_2\ldots i_m}>0$ for some $1 \le \alpha \le \beta$. Hence, $\ds\sum_{\alpha=1}^\beta {\cal A}^\alpha > 0$ and, subsequently, 
${\cal B}^\beta > 0$. 
\ep

As an illustration of Theorem \ref{irr2pri}, let us look at the following:

\begin{6}
\label{irrpriex2}
We consider the third order, $4$ dimensional tensor below, which is obtained by adding $\cal I$ of the same size to the s-irreducible tensor in Examples \ref{irrpriex0} and \ref{irrpriex1}.
$${\cal A}(:,:,1)=\left[\begin{array}{cccc}
2 & 0 & 0 & 0\\
1 & 1 & 1 & 0\\
0 & 1 & 1 & 1\\
0 & 0 & 0 & 1
\end{array}\right], ~{\cal A}(:,:,2)=\left[\begin{array}{cccc}
1 & 0 & 1 & 1\\
0 & 2 & 0 & 0\\
0 & 1 & 1 & 0\\
1 & 0 & 1 & 1
\end{array}\right],$$
$${\cal A}(:,:,3)=\left[\begin{array}{cccc}
1 & 1 & 0 & 1\\
1 & 1 & 1 & 0\\
0 & 0 & 2 & 0\\
0 & 0 & 0 & 1
\end{array}\right], ~{\cal A}(:,:,4)=\left[\begin{array}{cccc}
1 & 0 & 0 & 0\\
1 & 2 & 1 & 0\\
0 & 0 & 1 & 1\\
0 & 0 & 0 & 2
\end{array}\right].$$
Then, $\cal A$ is s-primitive. In fact, it can be easily verified that ${\cal A}^3 > 0$.
\end{6}

In the same vein as Theorem \ref{irr2pri}, we can also state: 
\begin{2}
\label{irr2pri2}
For a nonnegative $m$th order, $n$ dimensional tensor, if it is s-irreducible and if all of its diagonal entries are positive, then it must be s-primitive.
\end{2}
\bp
Denote such a tensor by $\cal A$. Let $c$ be the smallest diagonal entry of $\cal A$, i.e., 
$$c=\min_{i, i_3, \ldots, i_m \in \langle n \rangle} a_{iii_3\ldots i_m} > 0.$$
Obviously, we have 
$${\cal A} \ge \frac{1}{2}{\cal A}_d+\frac{1}{2}{\cal A} \ge \frac{c}{2}{\cal I}+\frac{1}{2}{\cal A}.$$
According to Lemma \ref{comp} and Theorem \ref{irr2pri}, therefore, the conclusion is now obvious.
\ep

To end this section, we recall that in the context of higher order Markov chains, s-irreducibility and s-primitivity translate into ergodicity and regularity, respectively. These concepts are essential in the study of such chains. Let $X$ be an $(m-1)$th order, $n$-state Markov chain. If $X$ is ergodic, then its mean first passage time tensor $\mu$, for instance, is well defined and is uniquely determined via \cite{HX24a, Xu26a}
$$\mu={\cal E}+(\mu -\mu_d)\boxtimes {\cal P},$$
where $\cal P$ is the transition tensor and $\cal E$ is the tensor of all ones of the same size as $\cal P$. On the other hand, the regularity of $X$ plays an important role in studying the limiting probability distribution of $X$. These results, however, may fail without their respective conditions of ergodicity or regularity. Theorems \ref{irrpri}, \ref{irr2pri}, and \ref{irr2pri2}, therefore, have practical implications as well. When $X$ is already known to be ergodic, for example, a sufficient condition as formulated in Theorem \ref{irr2pri2} for $X$ to be regular is just the positivity of the diagonal entries of $\cal P$.

\section{Conclusions}
\label{concl}
\setcounter{equation}{0}

In this work, we have established an alternative framework for irreducibility and primitivity of nonnegative tensors. Unlike the existing notions built on the zero-nonzero pattern of entries or the codomain of a nonlinear map, our framework is based entirely on the $\boxtimes$ product and power. While expanding the classical methodology on nonnegative matrices, such a framework opens up the door to a new way for further investigating irreducibility and primitivity of nonnegative tensors, distinct from existing works that have mostly focused on $Z$-eigenvalues and $Z$-eigenvectors. Seeing a vast body of classical results on nonnegative matrices, many intriguing questions are waiting to be explored. Given a nonnegative $m$th order, $n$ dimensional s-primitive tensor $\cal A$, for instance, we may ask what a potential upper bound on $\alpha$ is so as to guarantee ${\cal A}^\alpha > 0$. On the other hand, if instead such an $\cal A$ is known only to be s-irreducible, we may also ask whether s-primitivity can be archieved when some, but not all, of the zero diagonal entries of $\cal A$ are made positive. In addition, exploring additional applications of the results derived from our framework is another worthwhile direction for future research.

\end{document}